\documentclass[11pt,english]{article}
\usepackage[T1]{fontenc}
\usepackage[latin9]{inputenc}
\usepackage{amsmath}
\usepackage{amsthm}
\usepackage{amssymb}
\usepackage{graphicx}

\makeatletter
  \theoremstyle{plain}
  \newtheorem*{prop*}{\protect\propositionname}
  \theoremstyle{plain}
  \newtheorem*{thm*}{\protect\theoremname}
\theoremstyle{plain}
\newtheorem{thm}{\protect\theoremname}
  \theoremstyle{plain}
  \newtheorem{prop}[thm]{\protect\propositionname}
  \theoremstyle{definition}
  \newtheorem{example}[thm]{\protect\examplename}
  \theoremstyle{remark}
  \newtheorem{rem}[thm]{\protect\remarkname}
  \theoremstyle{plain}
  \newtheorem{cor}[thm]{\protect\corollaryname}

\usepackage{mathrsfs}\usepackage{amsthm}\usepackage{amscd}

\usepackage{hyperref}
\usepackage{graphics}
\usepackage{a4wide}
\@ifundefined{definecolor}
 {\usepackage{color}}{}

\newcounter{EQNR}
\renewcommand{\contentsname}{Contents\\
{\footnotesize\normalfont(A table of contents should normally not
be included)}}

\makeatother

\usepackage{babel}
  \providecommand{\corollaryname}{Corollary}
  \providecommand{\examplename}{Example}
  \providecommand{\propositionname}{Proposition}
  \providecommand{\remarkname}{Remark}
  \providecommand{\theoremname}{Theorem}
\providecommand{\theoremname}{Theorem}

\begin{document}

\title{Hahn-Banach for metric functionals and horofunctions}

\author{\date{November, 2019}Anders Karlsson \thanks{Supported in part by the Swiss NSF.}}
\maketitle
\begin{abstract}
It is observed that a natural analog of the Hahn-Banach theorem is
valid for metric functionals but fails for horofunctions. Several
statements of the existence of invariant metric functionals for individual
isometries and 1-Lipschitz maps are proved. Various other definitions,
examples and facts are pointed out related to this topic. In particular
it is shown that the metric (horofunction) boundary of every infinite
Cayley graphs contains at least two points. 
\end{abstract}

\section{Introduction}

It is well-recognized that the Hahn-Banach theorem concerning extensions
of continuous linear functionals is a cornerstone of functional analysis.
Its origins can be traced to so-called moment problems, and in addition
to H. Hahn and Banach, Helly's name should be mentioned in this context
(\cite{Di81,P07}).

In an important recent development, called the \emph{Ribe program,
}notions from geometric Banach space theory are formulated purely
in terms of the metric associated with the norm, and then studied
for metric spaces. This has been developed by Bourgain, Ball, Naor
and others, for a recent partial survey, see \cite{N18}. This subject
is in part motivated by significant applications in computer science.

It has been remarked in a few places (see \cite{Y11,Ka02} for merely
two references, the earliest discussion would surely be Busemann's
work in the 1930s, see \cite{Pa05} for a good exposition of Busemann's
work in this context) that Busemann functions or horofunctions serve
as replacement for linear functions when the space is not linear.
This has been a useful notion in the theory of manifolds with non-negative
curvature as well as for spaces of non-positive curvature, and more
recently several other contexts beyond any curvature restriction,
for example \cite{LN12,LRW18,W18,W19,KaL07}. In parallel, they have
been identified or described for more and more metric spaces, see
\cite{Ka19} for references.

In the present paper, we continue to consider the analogy in this
direction between the linear theory and metric theory. More precisely,
we consider the category of metric spaces and semi-contractions (non-expanding
maps or $1$-Lipschitz maps), that is, maps $f:X\rightarrow Y$ which
do not increase distances:
\[
d(f(x_{1}),f(x_{2}))\leq d(x_{1},x_{2})
\]
for all $x_{1},x_{2}\in X$. As discussed below one has metric analogs
of the norm, spectral radius, weak topology, the Banach-Alaoglu theorem,
and the spectral theorem. The present paper focuses on pointing out
an analog of the Hahn-Banach theorem. To our mind, the circumstance
that fruitful such analogs are present is surprising and promising. 

We follow Banach in his paper \cite{Ba25} and in his classic text
\emph{The Theory of Linear Operators} from 1932 in calling maps from
a metric space $X$ into $\mathbb{R}$ \emph{functionals. }Since we
in addition consider semi-contractions as our morphisms, we prefix
the functionals that we will consider by the word \emph{metric}. These
\emph{metric functionals} generalize Busemann functions and horofunctions,
see the next section for precise definitions of the latter two concepts,
and are moreover an analog and replacement for continuous linear functionals
in standard functional analysis, perhaps further motivating the use
of the word \emph{functional}.

To be precise, with a fixed origin $x_{0}\in X$, the following functions
and their limits in the topology of pointwise convergence are called
\emph{metric functionals :
\[
h_{x}(\cdot):=d(\cdot,x)-d(x_{0},x).
\]
}

The closure of this continuous injection of $X$ into functionals,
is called the \emph{metric compactification }of $X$ and is compact
by the analog of Banach-Alaoglu, see below. If $X$ is proper and
geodesics we call the subspace of new points obtained with the closure
the \emph{metric boundary. }In discussions with Tobias Hartnick, we
observed the following statement (and as will be explained the corresponding
statement for horofunctions is false). 
\begin{prop*}
(Metric Hahn-Banach statement) Let $(X,d)$ be a metric space with
base point $x_{0}$ and $Y$a subset containing $x_{0}.$ Then for
every metric functional $h$ of $Y$ there is a metric functional
$H$ of $X$ which extends $h$ in the sense that $H|_{Y}=h$. 
\end{prop*}
Because the above is a direct consequence of compactness, it is a
bit artificial to speak of applications of it in a strict sense. However,
it does provide a point of view that we believe is useful. 

Note that general Cayley graphs are highly non-trivial metric spaces,
so the following statement is not an obvious fact:
\begin{thm*}
Let $\Gamma$ be a finitely generated, infinite group and $X$ a Cayley
graph defined by $\Gamma$ and a finite set of generators. Then the
metric boundary of $X$ contains at least two points, and each metric
functional is unbounded. 
\end{thm*}
Fixed point theorems are of fundamental importance in analysis. The
following results are related to fixed point statements:
\begin{thm*}
For any monotone distorted isometry $g$ of a metric space $X$, there
exists a metric functional that vanishes on the whole orbit $\left\{ g^{n}x_{0}:n\in\mathbb{Z}\right\} $.
\end{thm*}
Thus this result shows that the classical picture of parabolic isometries
preserving horospheres in hyperbolic geometry extends, see Figure
1. These two theorems are proved in section \ref{sec:Distortion}
where also some further results about metric functionals of general
Cayley graphs are discussed.
\begin{figure}

\caption{A classical parabolic isometry.\protect\includegraphics[scale=0.4]{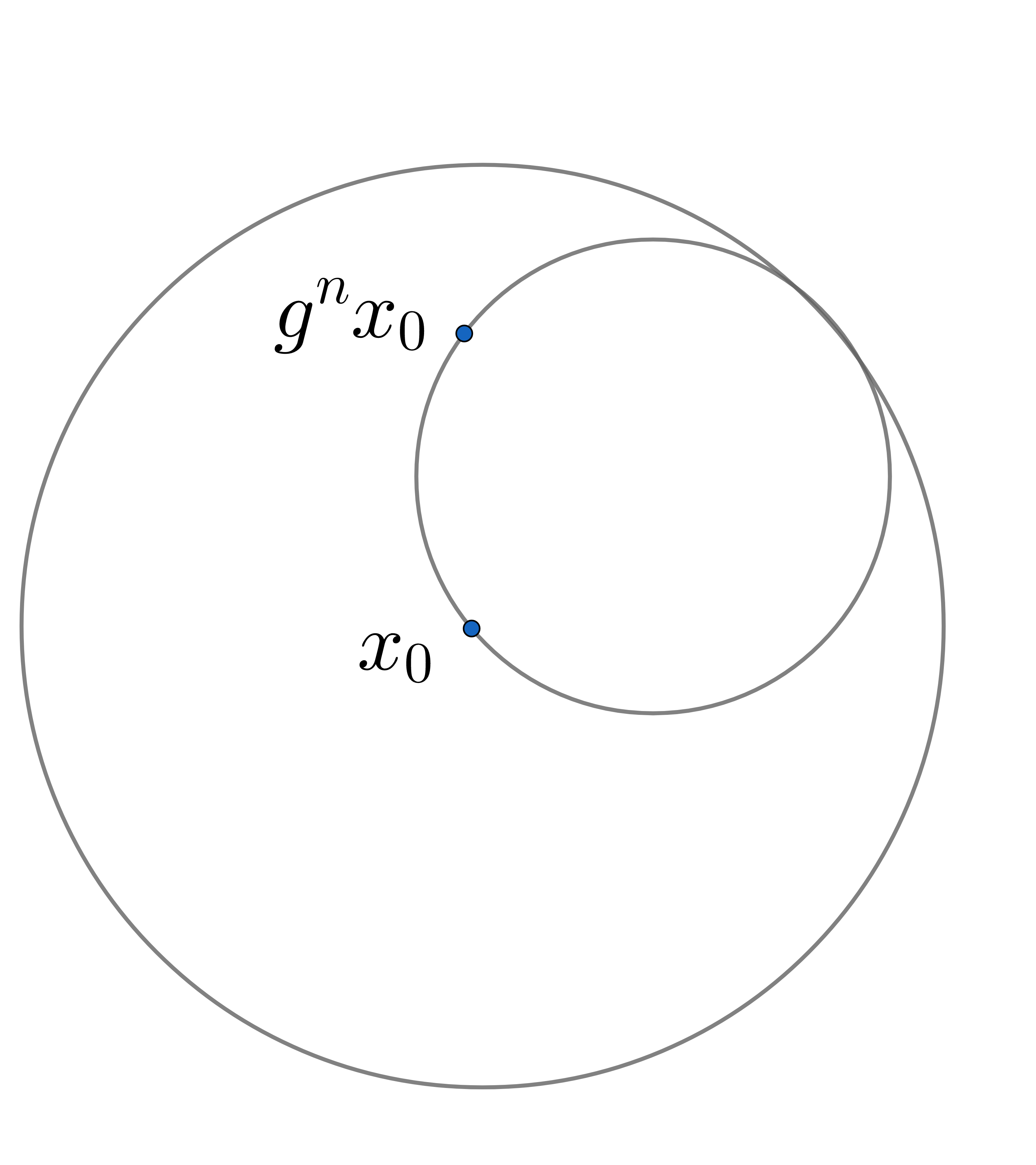}}

\end{figure}

With Bas Lemmens we observed the following improvement of the theorem
but under another hypothesis on the map:
\begin{prop*}
Suppose that $f$ is a semi-contraction of a metric space $X$ with
$\inf_{x}d(x,f(x))=0$. Then there is a metric functional $h,$ such
that 
\[
h(f(x))\leq h(x)
\]
 for all $x\in X$. In case $f$ is an isometry the inequality is
an equality.
\end{prop*}
It is an improvement in the sense that if $f$ is an isometry (with
the alternative assumptions) then $h(f^{n}x_{0})=h(x_{0})=0$ for
all $n$. 

Bader and Finkelshtein defined a \emph{reduced boundary }by identifying
functions in the metric boundary if they differ by a bounded function
\cite{BF19}. We note that:
\begin{prop*}
Every isometry of a metric space fixes a point in the reduced metric
compactification.
\end{prop*}

Busemann wrote in 1955 that ``... two startling facts: much of Riemannian
geometry is not truly Riemannian and much of differential geometry
requires no derivatives''. Although it would be too much to affirm
that much of functional analysis requires no linear structure, at
least there are a number of analogs for general metric spaces.

\textbf{Acknowledgements:} It is a pleasure to thank Tobias Hartnick
for the invitation to Giessen and our discussion there that sparked
the origin of this note. I also thank Nate Fisher, Pierre de la Harpe,
and Bas Lemmens for useful discussions.

\section{Basic definitions and terminology}

Let $(X,d)$ be a metric space. One defines (a variant of the map
considered by Kuratowski and Kunugui in the 1930s, see for example
\cite[p. 45]{L93}) 
\[
\Phi:X\rightarrow\mathbb{R\mathrm{^{X}}}
\]
via
\[
x\mapsto h_{x}(\cdot):=d(\cdot,x)-d(x_{0},x).
\]
As the notation indicates the topology we take here in the target
space is pointwise convergence. The map is continuous and injective. 
\begin{prop}
(Metric Banach-Alaoglu) The space $\overline{\Phi(X)}$ is a compact
Hausdorff space.
\end{prop}

\begin{proof}
By the triangle inequality we note that 
\[
-d(x_{0},y)\leq h_{x}(y)\leq d(y,x_{0}),
\]
which implies that the closure $\overline{\Phi(X)}$ is compact by
the Tychonoff theorem. It is Hausdorff, since it is a subspace of
a product space of Hausdorff spaces, indeed metric spaces. 
\end{proof}
In general this is not a compactification of $X$ in the strict and
standard sense that the space is embedded, but it is convenient to
still call it a compactification. In other words, we provide a weak
topology that has compactness properties.

By the triangle inequality it follows that all the elements in $\overline{X}:=\overline{\Phi(X)}$
are semi-contractive functionals $X\rightarrow\mathbb{R}$, and we
call them \emph{metric} \emph{functionals}.
\begin{example}
Let $\gamma$ be a geodesic line (or just a ray $\gamma:\mathbb{R}_{+}\rightarrow X$),
which is a standard notion in metric geometry at least since Menger.
Then the following limit exists:
\[
h_{\gamma}(y)=\lim_{t\rightarrow\infty}d(y,\gamma(t))-d(\gamma(0),\gamma(t)).
\]
\end{example}

The reason for the existence of the limit for each $y$ is that the
sequence in question is bounded from below and monotonically decreasing
(thanks to the triangle inequality), see \cite{BGS85,BrH99}. This
element in $\overline{\Phi(X)}$ is called the \emph{Busemann function
associated with }$\gamma$.
\begin{example}
The open unit disk of the complex plane admits the Poincaré metric,
which in its infinitesimal form is given by
\[
ds=\frac{2\left|dz\right|}{1-\left|z\right|^{2}}.
\]
This gives a model for the hyperbolic plane and moreover every holomorphic
self-map of the disk is a semi-contraction in this metric (the Schwarz-Pick
lemma). The Busemann function associated to the (geodesic) ray from
0 to the boundary point $\zeta$, in other words $\zeta\in\mathbb{C}$
with $\left|\zeta\right|=1$, is
\[
h_{\zeta}(z)=\log\frac{\left|\zeta-z\right|^{2}}{1-\left|z\right|^{2}}.
\]
These functions appear (in disguise) already in 19th century mathematics,
such as in the Poisson integral representation formula and in Eisenstein
series. 
\end{example}

The more common choice of topology, introduced by Gromov in \cite{Gr81},
is uniform convergence on bounded sets. We denote the corresponding
closure $\overline{X}^{h}$ and call it the \emph{horofunction bordification},
the new points are called \emph{horofunctions }following common terminology.
This closure amounts to the same compactification if $X$ is proper
(i.e. closed bounded sets are compact), see \cite{BrH99}, but in
general it is quite different, in particular there is no notion of
weak compactness. As in the example described in Remark \ref{rem:HBhorofunctions}
below, a space may have no horofunctions. Moreover, the following
example shows that Busemann functions are not always horofunctions,
since the limit above might not converge in this topology:
\begin{example}
\label{exa:Busemannray}Take one ray $\left[0,\infty\right]$ that
will be geodesic $\gamma$, then add an infinite number of points
at distance $1$ to the point $x_{0}=0$ and distance $2$ to each
other. Then at each point $n$ on the ray, connect it to one of the
points around 0 (that has not already been connected) with a geodesic
segment of length $n-1/2$. This way $h_{\gamma}(y)=\lim_{t\rightarrow\infty}d(\gamma(t),y)-d(\gamma(t),\gamma(0))$
still of course converge for each $y$ but not uniformly. Hence the
Busemann function $h_{\gamma}$ is a metric functional but not a horofunction.
See Figure 2.

\begin{figure}

\caption{The ray $\gamma$ does not define a horofunction.\protect\includegraphics[scale=0.4]{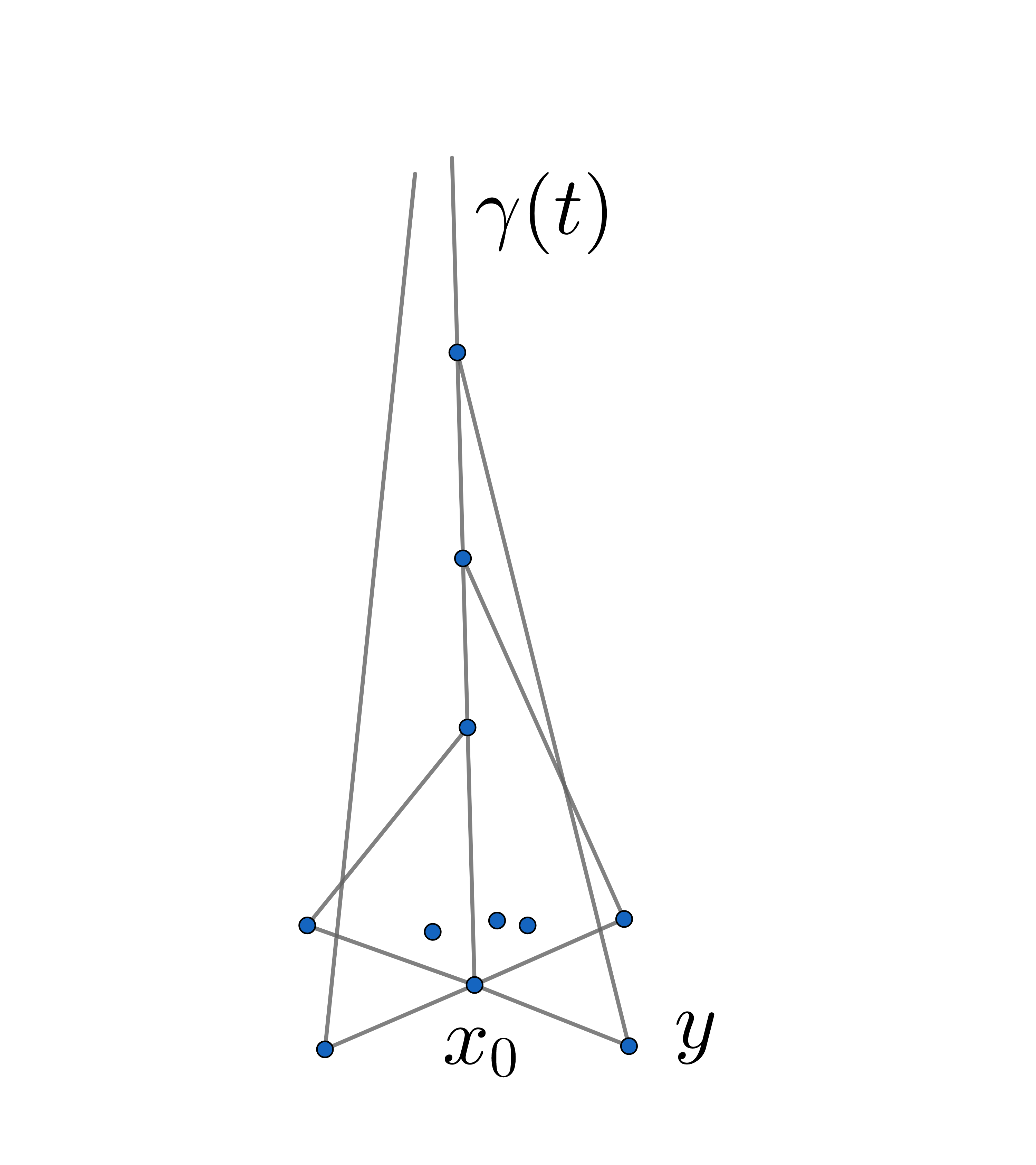}}
\end{figure}

The topology that we chose here has been useful in a few instances
already: \cite{GV12,GK15,G18,MT18} and in recent work by Bader and
Furman. 

In a Banach space the metric functional associated to points are not
linear since they are closely related to the norm, but sometimes their
limits are linear, see for example \cite{Ka19}. In any case they
are all convex functions:
\end{example}

\begin{prop}
Let $E$ be a normed vector space. Every function $h\in\overline{E}$
is convex, that is, for any $x,y\in X$ one has
\[
h(\frac{x+y}{2})\leq\frac{1}{2}h(x)+\frac{1}{2}h(y).
\]
\end{prop}

\begin{proof}
Note that for $z\in E$ one has
\[
h_{z}((x+y)/2)=\left\Vert (x+y)/2-z\right\Vert -\left\Vert z\right\Vert =\frac{1}{2}\left\Vert x-z+y-z\right\Vert -\left\Vert z\right\Vert 
\]

\[
\leq\frac{1}{2}\left\Vert x-z\right\Vert +\frac{1}{2}\left\Vert y-z\right\Vert -\left\Vert z\right\Vert =\frac{1}{2}h_{z}(x)+\frac{1}{2}h_{z}(y).
\]
This inequality passes to any limit point of such $h_{z}$. 
\end{proof}
In \cite[Lemma 3.1]{GK15} it is shown that for any metric functional
$h$ of a real Banach space there is a linear functional $f$ of norm
at most 1 such that $f\leq h$. The proof uses the standard Banach-Alaoglu
and Hahn-Banach theorems.
\begin{example}
Gutiérrez has provided a good description of metric functionals for
$L^{p}$ spaces $p\geq1$. To give an idea we recall the formulas
for $\ell^{p}(J)$ for $1<p<\infty$. There are two types:
\[
h_{z,c}(x)=\left(\left\Vert x-z\right\Vert _{p}^{p}+c^{p}-\left\Vert z\right\Vert _{p}^{p}\right)^{1/p}-c,
\]
where $z\in\ell^{p}(J)$ and $c\geq\left\Vert z\right\Vert _{p}$,
as well as 
\[
h_{\mu}(x)=-\sum_{j\in J}\mu(j)x(j),
\]
where $\mu\in\ell^{q}(J)$, with $q$ the conjugate exponent, and
$\left\Vert \mu\right\Vert _{q}\leq1$. See \cite{Gu17,Gu18,Gu19}
for more details and precise statements. An interesting detail that
Gutiérrez showed is that the function identically equal to zero is
not a metric functional for $\ell^{1}$ (in contrast to the space
$\ell^{2}$). He also observed how a famous fixed point free example
of Alspach fixes a metric functional. For another discussion about
Busemann functions of certain normed spaces, see \cite{W07,W18}. 
\end{example}

\section{Hahn-Banach theorem for metric functionals}
\begin{prop}
(Metric Hahn-Banach statement.) Let $(X,d)$ be a metric space with
base point $x_{0}$ and $Y$a subset containing $x_{0}.$ Then for
every $h\in\overline{Y}$ there is a metric functional $H\in\overline{X}$
which extends $h$ in the sense that $H|_{Y}=h$. 
\end{prop}

\begin{proof}
Given $h\in\overline{Y}$. Since the metric compactifications are
Hausdorff (even metrizable if $X$ is separable) we take a net $h_{y_{\alpha}}$
that converges to the unique limit $h.$ These points $y_{\alpha}$
are also points in $X$ and by compactness of $\overline{X}$ also
has a limit point $H$ there. By uniqueness of the limits it must
coincide with $h$. 
\end{proof}
\begin{rem}
\label{rem:HBhorofunctions}On the other hand, for horofunctions,
i.e. for $\overline{X}^{h}$, the Hahn-Banach theorem does \emph{not}
hold. Example \ref{exa:Busemannray} shows this, with $Y$ taken to
be the (image of the) geodesic ray $\gamma$. Then $Y$ clearly has
a metric functional $b$ that is a Busemann function, however in $X$
any sequence of points going to infinity (i.e. along $Y)$ cannot
converge in $\overline{X}^{h}$, but to extend $b$ we would need
such. Another illustration of this phenomenon (a counter-example to
the proof but not the statement) is the following: Consider longer
and longer finite closed intervals $\left[0,n\right]$ all glued to
a point $x_{0}$. See Figure 3. This becomes a countable (metric)
tree which is unbounded but contains no infinite geodesic ray. Denote
by $x_{n}$ the other end points of each interval. The sequence $h_{x_{n}}$
do not converge uniformly on balls. On the other hand it does so in
$Y$ being just the set $\left\{ x_{n}\right\} _{n\geq0}$. 
\end{rem}

\begin{figure}
\caption{A counter-example.\protect\includegraphics[scale=0.4]{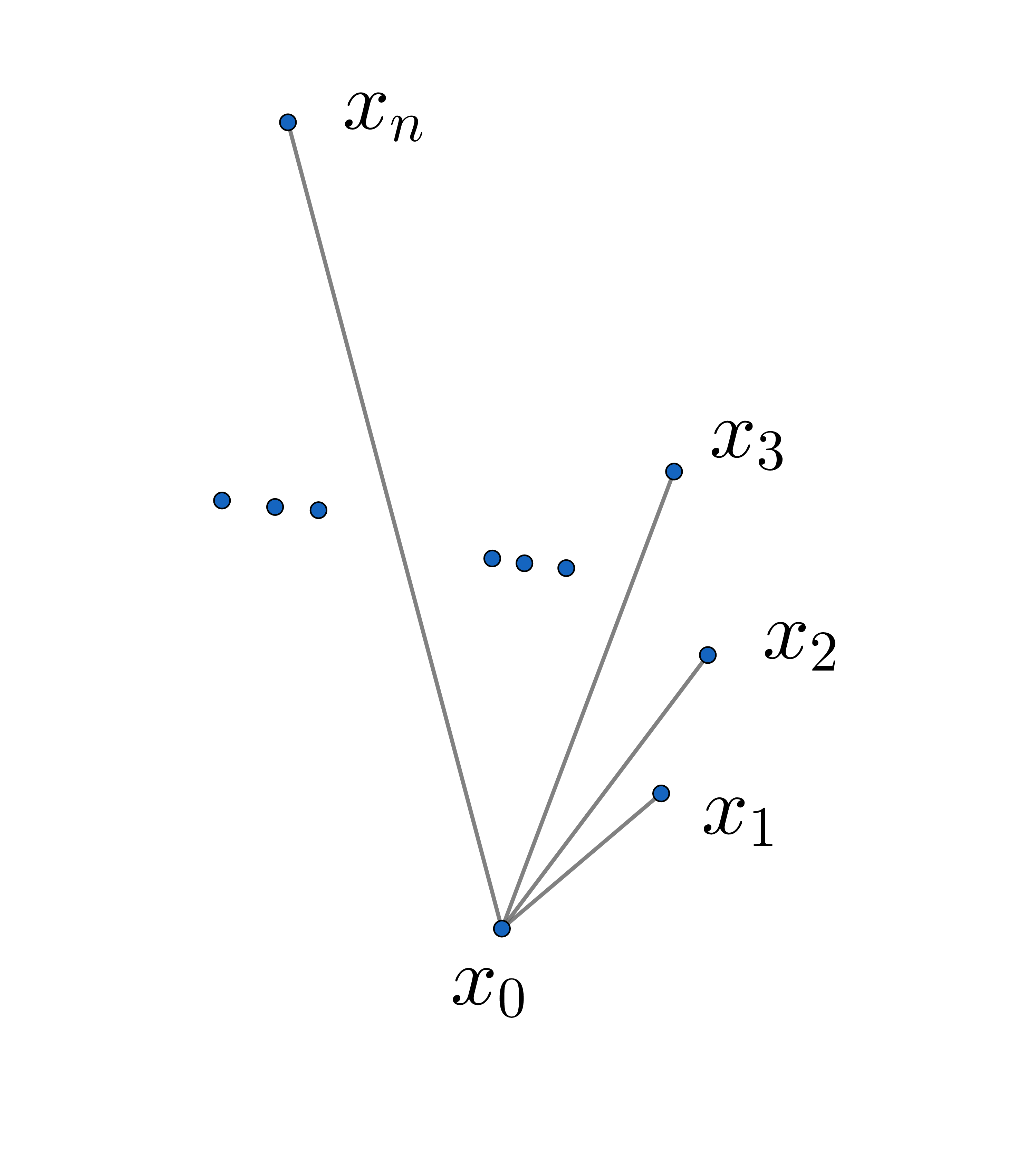}}
\end{figure}

Here is a remark in another direction: 
\begin{rem}
From a more general point of view, one could expect the possibility
of extensions of metric functionals. The real line $\mathbb{\left(R\mathrm{,}\left|\cdot\right|\right)}$
as a metric space is an injective object in the category we consider,
namely metric spaces and semi-contractions (i.e. $1$-Lipschitz maps),
in the following restricted sense: Given any metric space $A$ and
a monomorphism that is isometric $\phi:A\rightarrow B$ where $B$
is another metric space, and semi-contraction $f:A\rightarrow\mathbb{R}$
there is an extension in the obvious sense of $f$ to a morphism $B\rightarrow\mathbb{R}$,
for example
\[
\bar{f}(b):=\sup_{a\in A}\left(f(a)-d(\phi(a),b)\right)
\]
or 
\[
\bar{f}(b):=\inf_{a\in A}\left(f(a)+d(\phi(a),b)\right).
\]
To see this, in the former definition, note first that for $b'=\phi(a')$
one has $f(a')-d(\phi(a'),\phi(a'))\geq f(a)-d(\phi(a),\phi(a'))$,
which shows that our map is an extension, and then in general that
\[
\overline{f}(b)-\overline{f}(b')\leq\sup_{a}\left(f(a)-d(\phi(a),b)\right)-\sup_{a'}\left(f(a')-d(\phi(a'),b')\right)
\]
\[
\leq\sup_{a}\left(f(a)-d(\phi(a),b)-f(a)+d(\phi(a),b')\right)\leq d(b,b'),
\]
which implies it is a morphism. The origin of this observation is
\cite{McS34}.
\end{rem}

There are topological vector spaces with trivial dual. Maybe in spirit
this is a bit similar to the following example:
\begin{example}
\label{exa:R-distorted}Let $D:\mathbb{R}_{\geq0}\rightarrow\mathbb{R}_{\geq0}$
be an increasing function with $D(0)=0,$ $D(t)\rightarrow\infty$
and $D(t)/t\rightarrow0$ monotonically. The latter condition implies
that $D(t+s)\leq D(t)+D(s)$, so one sees that $(\mathbb{R},D(\left|\cdot\right|))$
is a metric space. As is observed in \cite{KaMo08}, $D(t)/t\rightarrow0$
implies that 
\[
\overline{(\mathbb{R},D(\left|\cdot\right|))}=\left\{ h_{x}:x\in\mathbb{R}\right\} \cup\left\{ h\equiv0\right\} .
\]
That is, there is a metric functional which vanishes on the whole
space, and the compactification is the one-point compactification.
Note that this metric space has no geodesics explaining in particular
why there are no Busemann functions.
\end{example}

In contrast to the linear theory, note that not every metric functional
of $X$ is a metric functional of a subset $Y$ as the following example
illustrates:
\begin{example}
Let $X$ be the Euclidean space $\mathbb{R}^{d}$ and $Y=\mathbb{R}$
be a one dimensional linear subspace. The Busemann function associated
to a ray from the origin perpendicular to $Y$ vanishes identically
on $Y$. On the other hand the zero function is not a metric functional
on $Y$. 
\end{example}

\section{Metric spectral notions}

Let $f$ be a semi-contraction of a metric space to itself. As remarked
in \cite{Ka19}, the \emph{minimal displacement} of $f$, $d(f)=\inf_{x}d(x,f(x))$,
is the analog of the norm of a linear operator and the \emph{translation
number} is the analog of the spectral radius:
\[
\ensuremath{\tau(f)=\lim_{n\rightarrow\infty}\frac{1}{n}d(x,f^{n}(x))},
\]
which exists in view of subadditivity. Similar to relationship between
the norm and spectral radius, one always has that $\tau(f)\leq d(f)$,
since the translation number is independent of $x$. Note that in
general the inequality may be strict, for example a rigid rotation
of the circle, or more interestingly, for groups with a word metric
all non-identity elements $g$have $d(g)\geq1$, but can easily have
$\tau(g)=0$. 

In passing we record the following simple fact.
\begin{prop}
The following tracial property holds: 
\[
\tau(gf)=\tau(fg)
\]
for any two semi-contractions $f$ and $g$. 
\end{prop}

\begin{proof}
From the triangle inequality, 
\[
d(x,(fg)^{n}(x))\leq d(x,f(x))+d(f(x),(fg)^{n}(x))
\]

\[
\leq d(x,f(x))+d(x,(gf)^{n-1}g(x))\leq d(x,f(x))+d(x,(gf)^{n}(x))+d((gf)^{n}(x),(gf)^{n-1}g(x))
\]

\[
\leq d(x,f(x))+d(x,(gf)^{n}(x))+d(f(x),x).
\]
Dividing by $n$ and sending $n$ to infinity shows one inequality.
By symmetry also the opposite inequality holds. 
\end{proof}
We say that an isometry $g$ is \emph{distorted }if $\tau(g)=0$.
Example \ref{exa:R-distorted} provides an observation related to
distortion. This example could be considered for $\mathbb{Z}$ (instead
of $\mathbb{R}$) and both are groups and the metric is invariant
(but not a word-metric). Here all non-zero elements are distorted.

It was shown in \cite{Ka01}, see also \cite{Ka19}, that for any
semi-contraction $f$ there is a metric functional $h$ such that
\[
h(f^{n}(x_{0}))\leq-\tau(f)n
\]
for all \emph{$n>0$} and 
\[
\lim_{n\rightarrow\infty}-\frac{1}{n}h(f^{n}(x))=\tau(f).
\]
I will refer to this result as the \emph{metric spectral} \emph{principle}.
As A. Valette pointed out to me, this could also be viewed as a statement
in the spirit of the classical Hahn-Banach theorem: existence of a
functional that realizes the norm of particular element. 

In discussions with Bas Lemmens we observed the following, which strengthen
the previous statement in a special case (another more general strengthening
appears in \cite{GV12}, but under weak non-positive curvature assumptions
on $X$):
\begin{prop}
Suppose that $f$ is a semi-contraction of $X$ with $d(f)=0$. Then
there is a metric functional $h,$ such that 
\[
h(f(x))\leq h(x)
\]
 for all $x\in X$. 
\end{prop}

\begin{proof}
For any $\epsilon>0$ we define the sets
\[
N_{\epsilon}=\left\{ x\in X:d(x,f(x))\leq\epsilon\right\} 
\]
and note that they are closed sets, and form a nested family. They
are all non-empty by the assumption on $f$. By compactness their
intersection must contain a metric functional $h$. There is thus
a sequence of points $x_{i}$ which converges to $h$ in the usual
weak sense (fix a sequence of $\epsilon$ going to $0$ and take a
point $x_{i}$ in each set). These $x_{i}$ are moved less and less
by $f$ which means that
\[
h(f(x))=\lim_{i\rightarrow\infty}d(f(x),x_{i})-d(x_{0},x_{i})=\lim_{i\rightarrow\infty}d(f(x),f(x_{i}))-d(x_{0},x_{i})
\]

\[
\leq\lim_{i\rightarrow\infty}d(x,x_{i})-d(x_{0},x_{i})=h(x).
\]
\end{proof}
\begin{cor}
Suppose the displacement of an isometry $g$ is zero, then there is
an invariant metric functional $h,$ in the strong sense that $h(gx)=h(x)$
for all $x$. 
\end{cor}

\begin{proof}
This is immediate from the previous proposition and its proof.
\end{proof}
Note that this is not in contradiction with a complicated isometry
like Edelstein's example, see \cite{Ka01}, since $h\equiv0$ is a
metric functional for Hilbert spaces. 

In the following section we will prove similar statements as above
in the case when $\tau(g)=0$ but $d(f)$ possibly strictly positive. 

\section{Group theory and parabolic isometries\label{sec:Distortion}}

Let me first record some facts that we, and probably others, have
realized years ago, see \cite{KaL07,Ka08}. They are however not generally
known by people in geometric group theory. Let $X$ be a Cayley graph
of a finitely generated group $\Gamma,$ which becomes a metric space
with the corresponding word metric associated to a finite generating
set, see the recent book \cite{CH16} for a wealth of metric geometry
in this setting. For Cayley graphs we always take the neutral element
as base point, $x_{0}=e$. The group acts by isometry on $X$, and
this actions extends continuously to an action by homeomorphism of
$\overline{X}$ and also on the \emph{metric boundary} $\partial X:=\overline{X}\setminus X.$
This boundary is a compact metrizable space. Let $\lambda$ be a $\Gamma$-invariant
probability measure on $\partial X$. Then the following map
\[
T(g):=\int_{\partial X}h(g)d\lambda(h)
\]
is a 1-Lipschitz group homomorphism $\Gamma\rightarrow\mathbb{R}$. 

No metric functional on $X$ is identically zero, see below. Therefore
if $\Gamma$ fixes a $h\in\partial X$ then there is a non-trivial
homomorphism $T:\Gamma\rightarrow\mathbb{Z}$. A further idea shown
in \cite{Ka08} is that given a finitely generated group with countable
boundary, there is a finite index subgroup that surjects on $\mathbb{Z}$
. One might wonder (\cite{Ka08}) or even conjecture (\cite{TY16})
that every group of polynomial growth has countable boundary. This
would immediately imply the celebrated theorem of Gromov that such
groups are virtually nilpotent. Some positive evidence for this approach
is provided in \cite{W11,TY16}. One could moreover have some hope
that growth considerations in relation to the metric boundary could
yield more than what is already known in this direction (recall that
Grigorchuk's Gap Conjecture remains open).

It is well-known that infinite finitely generated groups (and their
Cayley graphs) can be very complicated, for example it can contain
only elements of finite order, or even an exponent $N$ such that
$g^{N}=1$ for all group elements $g$. Therefore I think that the
following is a non-obvious fact.
\begin{thm}
Let $\Gamma$ be an infinite, finitely generated group and $X$ a
Cayley graph associated to a finite set of generators. Then the metric
boundary must contain at least two points and all metric functionals
are unbounded, in the sense that for any metric functional $h$ there
is no upper bound of $\left|h(x)\right|$ as a function of $x$. 
\end{thm}

\begin{proof}
First we establish the second assertion. Note that given the assumptions,
the graph $X$ is a proper metric space that is geodesic, in particular
connected, and with infinite diameter. In view of the latter property,
any $h_{g}$ with $g\in\Gamma$ is clearly unbounded. It remains to
consider a limit function $h(x)=\lim_{n\rightarrow\infty}d(x,g_{n})-d(e,g_{n})$.
Since the distance function is integer valued and the graph locally
finite, it holds that for any $r>0$ there is a number $N$ such that
$h(x)=d(x,g_{n})-d(e,g_{n})$ for all $n\geq N$ and $x$ such that
$d(x,e)\leq r.$ Take a geodesic from $e$ to $g_{N}$, it must intersect
the sphere around $e$ of radius $r$ in a point $y$. This means
that $h(y)=d(y,g_{N})-d(g_{N},e)=-r$. The statement follows since
$r$ was arbitrary. 

Now assume that the metric boundary contains exactly one point, $\partial X=\left\{ h\right\} .$
By the remarks above, $h$ defines a homomorphism $\Gamma\rightarrow\mathbb{Z}$
by 
\[
\gamma\mapsto h(\gamma).
\]
By what has just been shown, $h$ is unbounded, in particular not
identically $0$. This has as a consequence that the homomorphisms
is non-trivial. This means that there must exist an undistorted, infinite
order element $g$ (since any homomorphism into $\mathbb{Z}$ must
annihilate finite order and distorted elements), with $h(g)\neq0$.
Without loss of generality we may assume $h(g)>0$ (otherwise replacing
$g$by its inverse). But now by the metric spectral principle recalled
above, there must exist a metric functional $h_{1}$ such that $h_{1}(g)\leq0.$
This shows that $h_{1}$ is different from $h$. Since $\tau(g)>0$
the functional $h_{1}$ must take arbitrarily large negative values,
showing that $h_{1}\in\partial X$. 
\end{proof}
In classical hyperbolic geometry, parabolic isometries are those which
are distorted and preserves a horosphere. We say that an isometry
$g$ of a metric space $X$ with base point $x_{0}$ is \emph{monotone}
if $d(x_{0},g^{n}x_{0})\rightarrow\infty$ monotonically for all sufficiently
large $n$as $n\rightarrow\infty$. (So one way of treating a general
isometry $f$ might be to pass to a power of it, $g:=f^{N}$.) 
\begin{thm}
Let $g$ be a monotone isometry of a metric space $X$. If $\tau(g)=0$,
then there exists a metric functional that vanishes on the whole orbit
$\left\{ g^{n}x_{0}:n\in\mathbb{Z}\right\} $. 
\end{thm}

\begin{proof}
Let $Y$ be the orbit $\left\{ g^{n}x_{0}:n\in\mathbb{Z}\right\} $
with the metric induced by $X.$ It is a proper metric space since
$g$ is monotone. Consider the subset $A$ of $\partial Y$ of metric
functionals $h$ for which $h(g^{n}x_{0})\leq h(g^{m}x_{0})$ for
all $n\geq m.$ The subset $A$ is closed since the inequalities pass
to limits. It is also invariant under the group $H$ generated by
$g$since for any $n\geq m$,
\[
(g.h)(g^{n}x_{0})=h(g^{n-1}x_{0})-h(g^{-1}x_{0})=h(g^{n-1}x_{0})-h(g^{m-1}x_{0})+(g.h)(g^{m}x_{0})\leq(g.h)(g^{m}x_{0}).
\]

Next we verify that $A$ is non-empty. Take a converging subsequence
of $h_{g^{n}x_{0}}$ as $n\rightarrow\infty$, in notation $h_{g^{n_{i}}x_{0}}\rightarrow h$.
Then notice that by the monotonicity of $g$, for fixed $m\geq k$,
we have
\[
h(g^{m}x_{0})=\lim_{i}d(g^{m}x_{0},g^{n_{i}}x_{0})-d(x_{0},g^{n_{i}}x_{0})
\]
\[
=\lim_{i}d(x,g^{n_{i}-m}x_{0})-d(x_{0},g^{n_{i}}x_{0})
\]
\[
\leq\lim_{i}d(x_{0},g^{n_{i}-k}x_{0})-d(x_{0},g^{n_{i}}x_{0})=h(g^{k}x_{0}).
\]

Since $H$is a cyclic group acting on the compact non-empty set $A$
by homeomorphisms, there is an invariant probability measure $\mu$on
$A$. Therefore, as remarked above (with details found in \cite[Proposition 2]{Ka08}),
\[
T(g)=\int_{\partial Y}h(gx_{0})d\mu(h)
\]
defines a $1$-Lipschitz group homomorphism $T:H\rightarrow\mathbb{R}.$
Since any element of $\mathbb{R}$ is undistorted the image of $g$
must be $0,$ and so for every $n>0$ 
\[
0=\int_{\partial Y}h(g^{n}x_{0})d\mu(h).
\]

On the other hand $h(g^{n}x_{0})\leq h(x_{0})=0.$ This implies that
for every $n$ the set of $h$ for which $h(g^{n}x_{0})=0$ has full
measure. The intersection of countable full measure sets has full
measure, therefore there exists at least one $h$ which vanishes on
the whole orbit, thus proving the theorem.
\end{proof}
\begin{rem}
Note that the function $h\equiv0$ is a metric functional on any infinite
dimensional Hilbert space, while as an additive group no element is
distorted.
\end{rem}

Now we observe one thing in relation to the interesting notion of
reduced boundary from \cite{BF19} (note a conjecture in this paper
that states that for finitely generated nilpotent groups all reduced
boundary points should be fixed by the whole group). We extend their
definition by also considering points in $X$ and in the weak topology.
It is easy to verify that the action of isometries extends to the
reduced metric boundary, with the following calculation: Say that
$H$ is equivalent to $h$ differing by at most a constant $C$. Then
\[
\left|gH(x)-gh(x)\right|=\left|H(g^{-1}x)-H(g^{-1}x_{0})-h(g^{-1}x)+h(g^{-1}x_{0})\right|\leq2C.
\]

\begin{prop}
Let $g$ be an isometry of a metric space. Then it fixes a point in
the reduced metric compactification.
\end{prop}

\begin{proof}
Take any limit point $h$ of the orbit $g^{n}x_{0}$ which exists
by compactness. Recall that $(g.h)(x)=h(g^{-1}x)-h(g^{-1}x_{0}).$
The last term is bounded so we can forget this when passing to the
reduced compactification. We calculate:
\[
h(g^{-1}x)=\lim_{i\rightarrow\infty}d(g^{-1}x,g^{n_{i}}x_{0})-d(x_{0},g^{n_{i}}x_{0})
\]

\[
\leq\liminf_{i}d(g^{-1}x,g^{n_{i}-1}x_{0})+d(g^{n_{i}-1}x_{0},g^{n_{i}}x_{0})-d(x_{0},g^{n_{i}}x_{0})
\]

\[
=d(g^{-1}x_{0},x_{0})+\lim_{i}d(x,g^{n_{i}}x)-d(x_{0},g^{n_{i}}x_{0})=d(g^{-1}x_{0},x_{0})+h(x),
\]
for any $x\in X$. The reverse inequality is obtained by applying
the inequality to $gx$ instead of $x.$ Since the action by the isometry
$g$ is a well-defined map of the reduced boundary, this proves the
proposition.
\end{proof}

\section{Concluding remarks}

This short sections provide some remarks and suggestions of preliminary
nature. Recall that in Banach spaces weakly convergent sequences are
bounded. Also if a weakly convergent sequence belong to a closed convex
set, then the weak limit is also contained in this set (Mazur's theorem,
see for example \cite[p. 103]{La02}). Notice that in the example
in Remark \ref{rem:HBhorofunctions}, $h_{x_{n}}$ converges weakly
but is not bounded. On the other hand the limit belongs to the convex
hull of the sequence. One could wonder if there is a generalization
of this. In a yet different direction related to this, we refer the
reader to \cite{Mo16}. 

One could imagine defining the weak topology by declaring a sequence
\emph{$x_{n}$ }weakly converging to $x$ if $h(x_{n})\rightarrow h(x)$
for every metric functional on $X.$ But this would give back the
the usual strong, or metric, topology, since we could look at $h=h_{x}$.
Then it would make more sense to only consider $h$ which are at infinity.
(I am indebted to V. Guirardel and T. Hartnick for these remarks.)
Now it could be interesting to see in what way these $h$ coordinates
the space $x.$ for example if $X$ is a geodesic space moreover with
the property that every geodesic segment can be extended to a ray,
then one would have 
\[
d(x_{0},x)=\sup_{h}\left|h(x)\right|,
\]
where the supremum is taken over all Busemann functions. 

One could also try to define a norm on the Busemann functions (or
metric functionals coming from unbounded sequences):
\[
\left\Vert h\right\Vert =\sup_{x\neq x_{0}}\frac{\left|h(x)\right|}{d(x,x_{0})}.
\]
Note that this is not always $1,$ for example for infinite dimensional
Hilbert spaces where all linear functionals of norm at most 1 are
metric functionals, see \cite{Ka19}. With this one would have a metric
on the ``dual space'' of a metric space. At times, I have also suggested
the notion for a metric space to be reflexive: if every horofunction
is a Busemann function.

Here is a remark from \cite{Ka19}: In the works of Cheeger and collaborators
on differentiability of functions on metric spaces, the notion of
generalized linear function appears. In \cite{Ch99} it is connected
to Busemann functions, on the other hand that author remarks in \cite{Ch12}
that non-constant such functions do not exists for most spaces. Perhaps
it remains to investigate how metric functionals relate to this subject. 

It is shown above (and easily observed in case of abelian groups for
example) that the metric boundary always contains at least two points.
One may therefore conclude that is is not the Poisson boundary of
random walks. But it conceivable at least for some large classes of
groups, that if we pass to the reduced boundary, then random walks
could converge to points in this space (recurrent, or drift zero,
random walk could be said to converge to the interior which reduced
is one point). In this way hitting probability measures on the reduced
metric boundary could describe the behavior of random walks in a more
refined way. For example, random walks on the integers is governed
by the classical law of large numbers under a first moment condition.
The expectation value can be negative, zero or positive, depending
on the if the defining random walk measure is asymmetric. The reduced
metric compactification (with respect to any generating set) consists
of three points, that naturally can be denoted $-\infty,$ $+\infty$
and all finite points, that is, $\mathbb{Z}$ itself. And these three
points describe the asymptotic behavior of random walks with negative,
positive, or no drift, respectively. If the reduced boundary is not
enough one could consider the star geometry, in the sense of \cite{Ka05},
associated to the (reduced) metric compactification. This could be
invariant of the chosen generating set for the Cayley graphs. It would
also be of interest to study how groups act on this, their associated
incidence geometry. This is related to the conjecture in \cite{BF19}
already mentioned above.

Section de mathématiques, Université de Genève, 2-4 Rue du Lièvre,
Case Postale 64, 1211 Genève 4, Suisse 

e-mail: anders.karlsson@unige.ch 

and

Matematiska institutionen, Uppsala universitet, Box 256, 751 05 Uppsala,
Sweden 

e-mail: anders.karlsson@math.uu.se
\end{document}